\documentclass[11pt]{article}  
\usepackage{amssymb, amscd, amsmath, amsthm}

\newcommand{\beq}{\begin{equation}}
\newcommand{\eeq}{\end{equation}}

\newtheorem{corollary}{Corollary}
\newtheorem{theorem}{Theorem}

\allowdisplaybreaks

\title{A consequence of the classical Lindel\"of Hypothesis for the density of zeros of general Dirichlet series}

\author{by\\
J\'anos Pintz\thanks{Supported by the National Research Development and Innovation Office of Hungary, NKFIH, K133819, K147153 and Excellence 151341.}}

\date{Dedicated to the 50\textsuperscript{th} anniversary of the death of Paul Tur\'an and the 85\textsuperscript{th} birthday of G\'abor Hal\'asz}

\begin{document}
\maketitle

\setcounter{equation}{0}
\numberwithin{equation}{section}

\setcounter{equation}{0}
\numberwithin{equation}{section}



\begin{abstract}
Supposing the Lindel\"of hypothesis for the Riemann zeta function we prove that the density hypothesis holds for $\text{\rm Re }s > 3/4$ for a large class of Dirichlet series.
The first result of similar type was reached by G\'abor Hal\'asz and Paul Tur\'an in 1969.
\end{abstract}

\section{Introduction}
\label{sec:1}

The very recently announced proof of the Riemann Hypothesis (RH),
\beq
\label{eq:1.1}
\zeta(s) \neq 0 \ \ \text{ for }\ \ \sigma \neq 1/2 \qquad (s = \sigma + it)
\eeq
(verified by Lean according to the announcement) has clearly many highly important consequences.
Thus, the Lindel\"of hypothesis (LH)
\beq
\label{eq:1.2}
\zeta(s) \ll_\varepsilon (|t| + 2)^\varepsilon\ \ \text{ for any }\ \ \sigma \geq 1/2, \ \  \varepsilon > 0,
\eeq
and trivially the Density hypothesis (DH)
\beq
\label{eq:1.3}
N(\sigma, T) := \sum_{\varrho; \beta \geq \sigma} 1 = O_\varepsilon \bigl(T^{2(1 - \sigma) + \varepsilon}\bigr) \ \text{ for any }  \varepsilon > 0,
\eeq
where $\varrho = \beta + i\gamma$ denotes the non-trivial zeros of the zeta function.

As Tur\'an formulates (\cite{Tur1984}, p. 366) it came as a surprise that Ingham \cite{Ing1937} proved that $LH$ implies $DH$.
Three decades later Hal\'asz and Tur\'an showed \cite{HT1969} the groundbreaking result that $LH$ implies the much stronger
\beq
\label{eq:1.4}
N(\sigma, T) < T^\varepsilon  \ \text{ for any }  \sigma > 3/4 \ \text{ and } \ \varepsilon > 0.
\eeq

In the same work they showed that for a large class of Dirichlet series $f(s)$ (different from the one we will treat in the present work) the analogue of \eqref{eq:1.4} is true for the number $N_f(\sigma, T)$ if {\bf both} the zeta function and the function $f(s)$ satisfy the Lindel\"of hypothesis.

In a recent work \cite{Pin2026} we showed for a large class of $\mathcal F$ of general Dirichlet series that the DH is true for $\sigma > 9/10$ without any hypothesis (using a large value theorem of Heath-Brown \cite{HB1979} which used his twelfth power moment estimate \cite{HB1978} and the same idea of Hal\'asz which was crucial in proving \eqref{eq:1.4}).

Our present goal is to show that supposing LH only for the Riemann zeta function, and therefore the truth of the announced, Lean verified RH implies the density hypothesis for $\sigma > 3/4$ for a large class of Dirichlet series.

The conditions for the elements of the class $\mathcal F$ of Dirichlet series will be

(i) \ \ a mean value estimate on the line $\sigma = 1/2$;

(ii) \ Ramanujan's condition;

(iii) a very weak growth condition for $\sigma \geq 1/2$.

\section{Statement of the result. Notation}
\label{sec:2}

As outlined at the end of Section \ref{sec:1} we define the class $\mathcal F$ of Dirichlet series $F(s)$ as those satisfying for $\sigma \geq 1/2$
\beq
\label{eq:2.1}
F(s) = \sum_{n = 1} a_n n^{-s} \qquad (s = \sigma + it) \ \ a_n \in C
\eeq
which is regular in the halfplane $\sigma \geq 1/2$ up to a possible pole of finite order at $s = 0$,
\beq
\label{eq:2.2}
F(s) \ll_{\varepsilon} e^{|s|^\varepsilon} \ \text { for any } \ \varepsilon > 0 \ \text{ and } \ \sigma \geq 1/2.
\eeq
Further, we suppose that for $a_n$ and the coefficients $a^*(n)$ of the series $1/F(s)$ we have
\beq
\label{eq:2.3}
a_n \ll n^\varepsilon, \ \ \ a^*(n) \ll n^\varepsilon \ \text{ for any }\ \varepsilon > 0.
\eeq
Finally we suppose
\beq
\label{eq:2.4}
\int\limits_{-T}^T \left|F\left(\frac12 + it\right)\right|^2 dt \ll_\varepsilon T^{1 + \varepsilon}
\ \text{ for any } \ \varepsilon > 0, \ T \geq 1.
\eeq

In the following let $\varrho = \beta + i\gamma$ denote the zeros of $F(s)$ (we do not work further on with zeros of $\zeta(s)$).
Let $\varepsilon$ be an arbitrary small positive generic constant and
\beq
\label{eq:2.5}
N_F(\sigma, T) := \sum_{\substack{\varrho; F(\varrho) = 0\\
\beta \geq \sigma, |\gamma| \leq T}} 1.
\eeq

\begin{theorem}
\label{th:1}
If the classical LH is true for $\zeta(s)$ and $\sigma \geq \theta$ (denoted by LH$(\theta)$) then we have for any function $F \in \mathcal F$
\beq
\label{eq:2.6}
N_F(\sigma, T) = \mathcal O_\varepsilon \left(T^{2(1 - \sigma) + \varepsilon}\right) \ \text{ for } \
\sigma > (1 + \theta)/2,
\eeq
i.e.\ the DH is true for $F(s)$ in the halfplane $\sigma > (1 + \theta)/2$.
\end{theorem}

\begin{corollary}
\label{cor:1}
If LH is true for $\zeta(s)$ with $\theta = 7/8$ then \eqref{eq:2.6} is true for $F \in \mathcal F$.
Specifically for $\theta = 7/8$, the DH is true for $\sigma > 15/16$ and for $\theta = 1/2$ for $\sigma > 3/4$.
\end{corollary}

Thus, accepting the Lean verified proof of RH, \ DH is unconditionally true for $\sigma > 3/4$ for all elements $F \in \mathcal F$.

\section{Proof of Theorem 1}
\label{sec:3}

We will use Montgomery's zero detection method \cite{Mon1971} and a simple but ingenious idea of Hal\'asz.
In view of the growth condition \eqref{eq:2.2} by Jensen's formula we have
\beq
\label{eq:3.1}
N_F(\theta, T + 1) - N_F(\theta, T) \ll T^\varepsilon.
\eeq
Hence, it is sufficient to prove the estimate \eqref{eq:2.6} for any set $\mathcal R$ of $T^\delta$ separated zeros ($\delta > 0$ small, but fixed arbitrary)
$\varrho_k = \beta_k + i\gamma_k = 1 - \eta_k + i\gamma_k$ with the property
\beq
\label{eq:3.2}
\left|\gamma_k - \gamma_\ell\right| \geq T^\delta, \ \ \ |\gamma_k| \leq T, \ \ \beta_k \geq \sigma := 1 - \eta.
\eeq

Let $|\mathcal R| = K$ and with a fixed but arbitrary small $\delta > 0$
\beq
\label{eq:3.3}
X = T^\delta, \ \ \mathcal L = \log T, \ \ \
M_X(s) = \sum_{n \leq X} a_n^* n^{-s}.
\eeq

The usual Mellin transform gives
\beq
\label{eq:3.4}
e^{-1/T} + \sum_{n > X} b_n n^{-\varrho_k} e^{-n/T} = \frac1{2\pi i} \int\limits_{(2)} F(\varrho_k + s) M_X (\varrho_k + s) T^s \Gamma(s)ds
\eeq
where $b_n$ is the $n\textsuperscript{th}$ coefficient of $G(s):=F(s) M_X(s)$

We have for the zeros $\varrho_k = \beta_k + i\gamma_k$ of $F(s)$
\beq
\label{eq:3.5}
I(\varrho_k)\! := \frac1{2\pi i}\! \int\limits_{(2)}\!\! G(\varrho_k + s)T^s \Gamma(s)ds \! =\!
\frac1{2\pi i}\!\! \int\limits_{(1/2 - \beta_k)}\!\!\!\! G(\varrho_k + s) T^s \Gamma(s) ds\ + \ r(X, T, \varrho_k)
\eeq
where $r(X, T, \varrho_k)$ denotes the residue of the integrand at $s = 1 - \varrho_k$ (since the integrand is regular at $s = 0$).
If $m$ is the order of the pole of $F(s)$ at $s = 1$, then
\beq
\label{eq:3.6}
r(X, T, \varrho_k) \ll \sum_{0 \leq j \leq m} \left|\bigl(M_X(\varrho_k + s) T^s \Gamma(s)\bigr)^{(j)}_{s = 1 - \varrho_k}\right| \ll T^\delta T^{1 - \beta_k} e^{-T} = o(1).
\eeq

Due to \eqref{eq:3.6} and the exponential decay of the $\Gamma$-function we have
\beq
\label{eq:3.7}
\widetilde I(\varrho_k): = \frac1{2\pi i} \int\limits^{1/2 - \varrho_k + i\mathcal L^2}_{1/2 - \varrho_k - i\mathcal L^2} G(\varrho_k + s) T^s \Gamma(s)ds = 1 + \sum_{n > X} b_n n^{-\varrho_k} e^{-n/T} + o(1).
\eeq

Let us define $\mathcal R_1 \subset \mathcal R$ as the class of zeros where the above integral is $>1/3$.
The mean value estimate \eqref{eq:2.4} and the separation condition \eqref{eq:3.2} imply
\beq
\label{eq:3.8}
\int\limits^{\gamma_k + \mathcal L^2}_{\gamma_k - \mathcal L^2}|G(1/2 + it)|dt \gg T^{\beta_k - 1/2 + o(1)} \gg T^{1/2 - \eta},
\eeq
\beq
\label{eq:3.9}
\int\limits^{\gamma_k - \mathcal L^2 \gamma_k + \mathcal L^2}_{\gamma_k - \mathcal L^2}|G(1/2 + it)|^2 dt \gg T^{1 - 2\eta + o(1)} ,
\eeq
hence, by $G(s) \ll F(s)T^\delta$ we have
\beq
\label{eq:3.10}
\left|\mathcal R_1\right| T^{1 - 2\eta} \ll T^{1 + \delta + \varepsilon} \Longrightarrow \left|\mathcal R_1\right| \ll T^{2\eta + \delta + \varepsilon}.
\eeq

We now turn our attention to the remaining zeros $\varrho_k \in \mathcal R_2 := \mathcal R \setminus \mathcal R_1$ which satisfy
\beq
\label{eq:3.11}
\biggl|\sum_{X < n < T\mathcal L^2} b_n n^{-\varrho_k} e^{-n/T}\biggr| > 1/2
\eeq
since by $b_n = n^{o(1)}$ we can neglect all terms with $n \geq T\mathcal L^2$.
Let $|\mathcal R_2| = K_2$.
We have by the choice $X = T^\delta$ in \eqref{eq:3.3} trivially with $X_1 = \sqrt{X}$ by $b_n = n^{o(1)}$ from \eqref{eq:3.11}
\beq
\label{eq:3.12}
\Bigl|\sum_{X < n < T\mathcal L^2} b_n n^{-\varrho_k} \left(e^{-n/T} - e^{n/X_1}\right)\Bigr| > \frac13.
\eeq

Further on we will use the following lemma, valid under the assumption that
\beq
\label{eq:3.13}
\zeta(\sigma + it) \ll (|t| + 2)^\varepsilon \ \ \text{ for } \ \sigma \geq \theta ,
\eeq
denoted by $LH(\theta)$.
Let $\nu(n) = e^{-n/T} - e^{-n/X_1}$, $a > 0$ fixed.

\medskip
\noindent
{\bf Lemma.} {\it Under $LH(\theta)$ we have for $\alpha > \theta$, $T^a \leq |\gamma| \leq 2T$}
\beq
\label{eq:3.14}
S := \sum_{X < n < T\mathcal L^2} \frac{\nu(n)}{n^{\alpha + i\gamma}} \ll T^{\theta - \alpha + o(1)}.
\eeq

We emphasize that the Lemma refers for a modified partial sum of $\zeta(s)$ (not $F(s)$).

\begin{proof}
Similarly to \eqref{eq:3.4}--\eqref{eq:3.5} we have now
\begin{align}
\label{eq:3.15}
S &= \frac1{2\pi i} \int\limits_{(2)} \zeta(\alpha + i\gamma + s) \left(T^s - X_1^s\right) \Gamma(s) ds =\\
&= \frac1{2\pi i} \int\limits_{-\infty}^\infty \!\! \zeta(\theta + i\gamma + it)
\bigl(T^{it + \theta - \alpha} - X_1^{it + \theta - \alpha}\bigr) \Gamma(it + \theta - \alpha) dt
 + O\left(e^{-T^{a/3}}\right)
\nonumber
\end{align}
since the integrand is regular at $s = 0$ and the residum at $s = 1 - \alpha - i\gamma$ is $\ll T^{1 - \alpha}|\Gamma(1 - \alpha - i\gamma)| \ll e^{-T^a}$.
Hence, by the exponential decay of the $\Gamma$ function we have by $LH(\theta)$ and $\alpha > \theta$
\beq
\label{eq:3.16}
S \ll T^{\theta - \alpha + o(1)} . \qquad \hfill \qedhere
\eeq
\end{proof}

Now we are able to apply Hal\'asz's idea to the sum \eqref{eq:3.12}, i.e.\ to write now for every $\varrho_k = \varrho = 1 - \eta_k + i\gamma_k \in \mathcal R_2$ with the notation $\nu(n) = e^{-n/T} - e^{-n/X_1}$
\beq
\label{eq:3.17}
\biggl|\sum_{X < n < T\mathcal L^2} b_n \nu(n) n^{-\varrho_k}\biggr| = \varphi_k \sum_{X < n < T\mathcal L^2} b_n\nu(n) n^{-\varrho_k}
\eeq
where $|\varphi_k| =  1$.
Hal\'asz's idea is to square the RHS of the sums in \eqref{eq:3.17} for every $\varrho_k \in \mathcal R_2$
$\left(\left|\mathcal R_2\right| = K_2\right)$, interchange the order of summation over $k$ and $n$ and use the Cauchy--Schwarz inequality when $n$ runs in $I:= \left(X, T\mathcal L^2\right)$ with
\beq
\label{eq:3.18}
b_n \nu(n)n^{-\varrho_k} = b_n \left(\frac{\nu(n)}{n}\right)^{1/2} \cdot \left(\frac{\nu(n)}{n}\right)^{1/2} n^{\eta_k - i\gamma_k} := d_n\cdot e_n^{(k)}.
\eeq

Separating the diagonal terms (those with $k = \ell$) we deduce from $b_n = n^{o(1)}$ and the Lemma (following from $LH((\theta))$
\begin{align}
\label{eq:3.19}
|\mathcal R_2|^2 = K_2^2 &\ll \biggl(\sum_{k = 1}^{K_2} \varphi_k \sum_{n \in I} b_n \nu(n)n^{-\varrho_k}\biggr)^2 = \biggl(\sum_{n \in I} d_n \sum_{k = 1}^{K_2} \varphi_k e^{(k)}_n\biggr)^2\\
&\ll \biggl(\sum_{n \in I} \frac{|b_n|^2 \nu(n)}{n}\biggr)
\biggl(\sum_{k = 1}^{K_2}\sum_{\ell = 1}^{K_2} \varphi_k \overline{\varphi_\ell} \sum_{n \in I} \frac{\nu(n)}{n^{1 - \eta_k - \eta_\ell - i(\gamma_\ell - \gamma_k)}}\biggr)
\nonumber\\
&\ll T^{o(1)} \left\{ K_2(K_2 - 1) T^{\theta + 1 - 2\sigma + o(1)} + K_2 \mathcal L^c T^{2\eta}\right\}.
\nonumber
\end{align}
Since by $\sigma > (1 + \theta)/2$ the first term in the last sum is (much) smaller than the LHS of \eqref{eq:3.19} we have
\beq
\label{eq:3.20}
\left|\mathcal R_2\right| = K_2 \ll T^{2\eta + o(1)}. \qquad \hfill \qedhere
\eeq

\section{An application of Theorem 1}
\label{sec:4}

If the coefficients $a_n$ of $F(s)$ represent a completely multiplicative arithmetic function, i.e. $a_{mn} = a_m a_n$, then $a_1 = 1$ and $a^*(n) = \mu(n)a_n$.
Thus it is sufficient to suppose the Ramanujan condition for $F(s)$.
In such a way if we consider a Beurling number system with the following three properties

(i) \ \ the generating elements are positive integers;

(ii) \ the counting function of the integers satisfies
$$
N(x) = \kappa x + O(x^\varepsilon) \ \text{ for any } \ \varepsilon > 0;
$$

(iii) the Ramanujan condition holds for the coefficients of the Beurling zeta function $\mathcal B(s)$ of the system

\noindent
then we have

\begin{corollary}
\label{cor:2}
If LH is true for $\zeta(s)$ for $\sigma \geq \theta$ then if {\rm (i)--(iii)} hold then DH is true for $\mathcal B(s)$ for $\sigma > (1 + \theta)/2$.
In particular, if LH is true for $\theta = 1/2$ then DH is true for $\mathcal B(s)$ for $\sigma > 3/4$.
\end{corollary}

\begin{proof}
The mean value condition is satisfied in this case by Theorem 1.1 of Broucke and Hilberdink \cite{BH2024}.
\end{proof}

This improves a joint result of S. G. R\'ev\'esz and the author \cite{PR2025}.

A further application is for all elements of the Selberg class with degree at most $2$.
In this case the mean value estimate is always true as it follows from the approximate functional equation for $F(s)$ and from Montgomery--Vaughan's mean value theorem for Dirichlet polynomials.
The Ramanujan condition is satisfied for $F(s)$ and one can get around of it for $1/F(s)$ by the simple Lemma 1 of \cite{KP2003}.
So we have

\begin{corollary}
If LH is true for $\varrho(s)$ for $\sigma \geq \theta$ then the DH is true for elements $F(s)$ of the Selberg class of degree at most $2$.
In particular if LH is true for $\zeta(s)$ for $\theta = 1/2$ then DH is true for the mentioned functions $F(s)$ for $\sigma > 3/4$.
\end{corollary}

The above Corollary improves the recent result of the author \cite{Pin202?}
valid unconditionally for $\text{\rm Re }s \geq 0.9$.

{\small
\noindent
J\'anos Pintz\\
HUN-REN Alfr\'ed R\'enyi Institute of Mathematics \\
Budapest, Re\'altanoda u. 13--15\\
H-1053 Hungary\\
e-mail: pintz@renyi.hu}

\end{document}